\documentclass[12pt,oneside,english]{amsart}
\usepackage[latin1]{inputenc}
\usepackage{amssymb}

\makeatletter
\newtheorem{theorem}{Theorem}
\newtheorem{lemma}{Lemma}
\newtheorem{example}{Example}

\newtheorem{remark}{Remark}

\newtheorem{corollary}{Corollary}

\usepackage{babel}

\makeatother
\begin{document}
\baselineskip=17pt

\title[New proof of Coquet's
theorem]
{Two algorithms for evaluation of the Newman digit sum, and a new proof of Coquet's
theorem}

\author{Vladimir Shevelev}
\address{Department of Mathematics \\Ben-Gurion University of the
 Negev\\Beer-Sheva 84105, Israel. e-mail:shevelev@bgu.ac.il}


\begin{abstract}
We give two simple algorithms for the evaluation of difference between the numbers
of multiples of 3 with even and odd binary digit sums in interval $[0,x),$ and give
an elementary proof of Coquet's sharp estimates for it. 
\end{abstract}

\maketitle

\section{INTRODUCTION }

Let $p> 1,\;q> 1,\;m> 1$ be integers,  $\gcd(p, q-1)=1$. Let in base $q$

$$n=\sum^ v_{k=0}a_k  q^k,\;\;0\leq a_k< q.$$

Denote the digit sum of $n$ in base $q$ by $s_q(n):$
$$ s_q(n)=\sum^v_{k=0}a_k.$$
In 1968, Gelfond \cite{2} obtained the following general result
about the distribution of digit sums of integers: the number of integers $n < x$ satisfying $n\equiv l \pmod
{m}\;\;s_q(n)\equiv t \pmod {p}$ equals $\frac{x}{mp}+
O(x^\lambda), \;\;\lambda < 1$, where $\lambda$ does not depend on
$x,\;m,\;l$ and $t$.

     In particular, in the case of $p=q=2$, Gelfond found that $\lambda=\frac{\ln 3}{\ln
4}$.  Thus in this case, for $m=3,$ we have

\begin{equation}\label {1}
\sum_{n\leq x:\;\;n\equiv l\pmod 3,\;\;s(n)\equiv 0\pmod 2} 1=\frac x 6
+O(x^\lambda)
\end{equation}
and

\begin{equation}\label {2}
\sum_{n\leq x:\;\;n\equiv l\pmod 3,\;\;s(n)\equiv 1\pmod 2} 1=\frac x 6
+O(x^\lambda)
\end{equation}

with $s(n)=s_2(n),\;\;\lambda=\frac{\ln 3}{\ln 4}$.\newline
 Everywhere below we we consider this special case and suppose that $s(n)=s_2(n)$ and $\lambda=\frac{\ln 3}{\ln
4}.$

 For $x\in \mathbb{N}$ and $m\geq 2$, denote by $S_{m,\;l}(x)$ the sum

\begin{equation}\label {3}
S_{m,\;l}(x)=\sum_{0\leq n < x: \;\;n\equiv l\pmod m}(-1)^{s(n)}.
\end{equation}
Note that, in particular, $S_{3,\;0}(x)$ equals the difference between the numbers of multiples of 3 with even and odd binary digit sums (or multiples of 3 from sequences A001969 and A000069 in \cite{4}) in interval $[0,x).$
From (\ref{1}) and (\ref{2}) it follows that

\begin{equation}\label{4}
S_{3,\;l}(x)=O(x^\lambda).
\end{equation}

Leo Moser (cf. \cite{3}, Introduction) conjectured that always

\begin{equation}\label{5}
S_{3,0}(x)>0.
\end{equation}

 Newman \cite{3} proved this conjecture. Moreover, he obtained
the inequalities

\begin{equation}\label{6}
\frac{1}{20} < S_{3,0}(x)x^{-\lambda}< 5.
\end{equation}

In 1983, Coquet \cite{1} studied a very complicated continuous and nowhere differentiable
function $F(x)$ with period 1 for
which

\begin{equation}\label{7}
S_{3,0}(3x)=x^\lambda F\left(\frac{\ln x}{\ln
4}\right)+\frac{\eta(x)}{3},
\end{equation}

where

\begin{equation}\label{8}
\eta(x)=\begin{cases} 0,\;\; if \; x \;\; is\;\; even,\\
(-1)^{s(3x-3)}, \;\; if \;\; x \;\; is\;\; odd.\end{cases}
\end{equation}

He obtained that

\begin{equation}\label{9}
\limsup_{x\rightarrow
\infty,\;x\in\mathbb{N}}S_{3,0}(3x)x^{-\lambda}=\frac{55}{3}\left(\frac{3}{65}\right)^\lambda=1.601958421...,
\end{equation}

\begin{equation}\label{10}
\liminf_{x\rightarrow\infty,\;x\in\mathbb{N}}S_{3,0}(3x)x^{-\lambda}=\frac{2\sqrt{3}}{3}
=1.154700538...\;.
\end{equation}
\begin{remark}\label{r1}
J.-P. Allouche informed the author about a misprint in (\ref{8}). It should be
$$ \eta(x)=\begin{cases} 0,\;\; if \; x \;\; is\;\; even,\\
(-1)^{s(3x-1)}, \;\; if \;\; x \;\; is\;\; odd.\end{cases}$$
\newpage
It is important, since in the published case, using the periodicity of $F(x),$ it could be shown that (\ref{7}) and (\ref{8}) contradict one to another.
\end{remark}

     In this paper, using an elementary method, we give a new proof of Coquet's
result. Moreover, we prove that
for $x\geq 2,\;\;x\in \mathbb{N}$,

\begin{equation}\label{11}
\frac{2\sqrt{3}}{3}x^\lambda\leq
S_{3,0}(3x)\leq\frac{55}{3}\left(\frac {3}{65}\right)^\lambda
x^\lambda.
\end{equation}
Our method is based on a simple algorithm for the exact
evaluation of $S_{3,0}(x)$ described in Section 2. It allows to calculate also
$S_{3,1}(x),\;S_{3,2}(x)$ and $S_{3\cdot 2^k, r}(x),\;r\leq 3\cdot
2^k -1$ and to obtain for them the corresponding sharp estimates. In Section 6 we also give a fast algorithm for exact evaluation of the Newman digit sum $S_{3,0}(N)$ which is very suitable for sharp experiments.

\section{An algorithm for exact evaluation of $S_{3,0}$}
Let $y=2^{k_1}+2^{k_2}+\ldots +2^{k_r}$ be the binary expansion of
$y\in\mathbb{N}$.

Put
\begin{equation}\label{12}
t(y)=(-1)^{k_1}+(-1)^{k_2}+\ldots +(-1)^{k_r}.
\end{equation}

Evidently,
\begin{equation}\label{13}
t(y)\equiv y\pmod 3.
\end{equation}

Denote

\begin{equation}\label{14}
S_{3,0}([y, y+z))=S_{3,0}(y+z)-S_{3,0}(y).
\end{equation}

Consider an integer $m\in[1, k_r)$. Below we use some trivial
bijections of shifts such that the number of the integers divisible
by $3$ does not  change. According to (\ref{13}), we see that the
calculation of $S_{3,0}([y, y+2^m))$ reduces to the following cases:

a) $t(y)\equiv 0\pmod 6.$ Then

\begin{equation}\label{15}
S_{3,0}([y, y+2^m))=S_{3,0}(2^m);
\end{equation}

b) $t(y)\equiv 1\pmod 6.$ Then for any $n>\frac m 2$, $t(y)\equiv
t(2^{2n})\pmod 3$. Therefore, in the bijection of shift $[y,
y+2^m)\longleftrightarrow[2^{2n},\;2^{2n}+ 2^m)$ the integers
divisible by $3$ in $[y, y+2^m)$ correspond to integers in
$[2^{2n},\;2^{2n}+ 2^m)$ which also multiple of $3$.

Thus we have
\begin{equation}\label{16}
S_{3,0}([y,\;y+2^m))=S_{3,0}([2^{2n},\;2^{2n}+ 2^m));
\end{equation}
\newpage
c) $t(y)\equiv 2\pmod 6.$ Then for any
$k>\frac{m+2}{2},\;n>\frac{m+1}{2}$ we have $t(y) \equiv
t(2^{2k}+2^{2k-2})\equiv t(2^{2n-1})\pmod 3$ and thus, as above, we
find
$$S_{3,0}([y,\;y+2^m))=S_{3,0}([2^{2k}+2^{2k-2},\;2^{2k}+2^{2k-2}+2^m))=$$
\begin{equation}\label{17}
-S_{3,0}([2^{2n-1},\;2^{2n-1}+2^m));
\end{equation}
d) $t(y)\equiv 3\pmod 6$. Then, evidently,
\begin{equation}\label{18}
S_{3,0}([y,\;y+2^m))=-S_{3,0}(2^m);
\end{equation}
e) $t(y)\equiv 4\pmod 6$. For any $k>\frac{m+6}{2},\;n>\frac m 2,$
we have $$t(y)\equiv t(2^{2k}+2^{2k-2}+2^{2k-4}+2^{2k-6})\equiv
t(2^{2n})\pmod 3$$ and
$$ S_{3,0}([y,\;y+2^m))=$$ $$S_{3,0}([2^{2k}+2^{2k-2}+2^{2k-4}+2^{2k-6},\;
2^{2k}+2^{2k-2}+2^{2k-4}+2^{2k-6}+2^m))=$$
\begin{equation}\label{19}
-S_{3,0}([2^{2n},\;2^{2n}+2^m));
\end{equation}

f) $t(y)\equiv 5 \pmod 6$. In this case, for any $k>\frac{m+2}{2},\;n>\frac {m+1}{2},$ we have $t(y)\equiv
t(2^{2k}+2^{2k-2})\equiv t(2^{2n-1})\pmod 3.$ Therefore,
$$
S_{3,0}([y,\;y+2^m))=-S_{3,0}([2^{2k}+2^{2k-2},\;
2^{2k}+2^{2k-2}+2^m))=
$$
\begin{equation}\label{20}
=S_{3,0}([2^{2n-1},\;2^{2n-1}+2^m));
\end{equation}
Thus we conclude that, in order to calculate $S_{3,0}([y,\;y+2^m)),$ it is sufficient to find
$S_{3,0}(2^m)$ and $S_{3,0}([2^n,\;2^n +2^m)).\;n>m$. Below (Section 3) we shall prove the following theorem.
\begin{theorem}\label{t1}
\begin{equation}\label{21}
 S_{3,0}(2^m)=\begin{cases} 2\cdot3^{\frac m 2 -1}, \;\; if\;\;
m\;\; is \;\; even,\\ 3^{\frac{m-1}{2}},\;\; if\;\; m\; \;is \;\;
odd,\;\; m\geq 1;\end{cases}
\end{equation}

$$ S_{3,0}([2^n, 2^n+2^m))=$$
\begin{equation}\label{22}
\begin{cases} 3^{\frac m 2 -1}, \;\;
if\;\; m\;\; is \;\; even,\;\;1\leq m\leq n-1\\
3^{\frac{m-1}{2}},\;\; if\;\; m\; \;and \;\;n\;\; are\;\; odd,\;\;
1\leq m\leq n-2;\\0,\;\; if \;\;m \;\; is\;\; odd,\;\; n \;\; is\;\;
even,\;\; 1\leq m\leq n-1.\end{cases}
\end{equation}
\end{theorem}

In  addition to Theorem \ref{t1}, note that, for odd $N\geq 1,$

\begin{equation}\label{23}
S_{3,0}([N-1,N))=\begin{cases}(-1)^{s(N-1)},\;N\equiv 1\pmod
3\\0, \; otherwise \end{cases}
\end{equation}

Using formulas (\ref{15})-(\ref{23}), one can easily calculate
$S_{3,0}(x)$ for any $x\in\mathbb{N}$.
\newpage
\begin{example}\label{example 1}

Newman mentioned in \cite{3}about numerical studies by I.Barrodale
and R.MacLeod that beard out of the Moser`s conjecture up to $500
000$ and obtained that $S_{3,0}(500 000)$ is  around $17 000$. Let
us do an exact calculation. Using (\ref{15}), (\ref{26}),
(\ref{20})-(\ref{22}) we have

$$
S_{3,0}(500
000)=S_{3,0}(2^{18}+2^{17}+2^{16}+2^{15}+2^{13}+2^8+2^5)=
$$
$$
=S_{3,0}(2^{18})+S_{3,0}([2^{18}, 2^{18}+2^{17}))+S_{3,0}(2^{16})+
$$
$$
+S_{3,0}([2^{16}, 2^{16}+2^{13}))+S_{3,0}(2^{13})+S_{3,0}([2^{13},
2^{13}+2^8))+
$$
$$
+S_{3,0}(2^5)= 2\cdot 3^8 +0 + 2 \cdot 3^7 + 0 + 3^6 + 3^3 + 3^2 =
18261.
$$
\end{example}

\section{Proof of Theorem \ref{t1}}

Similar to (\ref{14}), for $S_{m,l}(x)$ (\ref{3}), we denote

\begin{equation}\label{24}
S_{m,l}([y, y+z))=S_{m,l}(y+z)-S_{m,l}(y).
\end{equation}

First of all, for any $x, y\in \mathbb{N}$,  it is easy to see that

\begin{equation}\label{25}
S_{6,0}([2x, 2y))=S_{3,0}([x, y)).
\end{equation}

Thus, since

\begin{equation}\label{26}
S_{3,0}([2x, 2y))=S_{6,0}([2x, 2y))+S_{6,3}([2x, 2y)),
\end{equation}

then we have

\begin{equation}\label{27}
S_{6,3}([2x, 2y))=S_{3,0}([2x, 2y))-S_{3,0}([x, y)).
\end{equation}

Furthermore, by (\ref{25}), evidently,

\begin{equation}\label{28}
S_{6,1}([2x, 2y))=-S_{3,0}([x, y))
\end{equation}
and
\begin{equation}\label{29}
S_{6,2}([2x, 2y))=-S_{6,3}([2x, 2y)).
\end{equation}
According to (\ref{27}), this means that
\begin{equation}\label{30}
S_{6,2}([2x, 2y))=S_{3,0}([x, y))-S_{3,0}([2x, 2y)).
\end{equation}

Furthermore, using (\ref{30}), we have
$$S_{3,1}([x, y))=S_{6,2}([2x, 2y))=$$
\begin{equation}\label{31}
S_{3,0}([x, y))- S_{3,0}([2x,
2y)).
\end{equation}
\newpage
On the other hand,

\begin{equation}\label{32}
S_{3,1}([x, y))=S_{6,1}([x,y))+S_{6,4}([x, y))
\end{equation}

and, consequently, by (\ref{28}), (\ref{31}) and (\ref{32}), we find

$$
S_{6,4}([2x, 2y))=S_{3,1}([2x,2y))-S_{6,1}([2x,
2y))=$$$$S_{3,1}([2x,2y))+S_{3,0}([x, y))=
$$
\begin{equation}\label{33}
=S_{3,0}([2x, 2y))-S_{3,0}([4x,4y))+S_{3,0}([x, y)).
\end{equation}
Therefore, we conclude that
\begin{equation}\label{34}
S_{3,2}([x, y))=S_{3,0}([x,y))+S_{3,0}([2x,2 y))-S_{3,0}([4x, 4y)).
\end{equation}

Furthermore, from (\ref{31}) and (\ref{34}) it follows that
$$S_{3,0}([x,y))+ S_{3,1}([x,y))+S_{3,2}([x,y))=$$
\begin{equation} \label{35}
3S_{3,0}([x,y))-S_{3,0}([4x,4y)).
\end{equation}
Now note that, the left hand side of (\ref{35}) equals to

$$S_{1,0}([x,y))=\sum_{x\leq n< y}(-1)^{s(n)}=$$
\begin{equation}\label{36}
\sum_{0\leq n<
y}(-1)^{s(n)}-\sum_{0\leq n< x}(-1)^{s(n)}.
\end{equation}

Let us show that, if $x$ is an even positive integer, then

\begin{equation}\label{37}
\sum_{0\leq n< x}(-1)^{s(n)}= 0.
\end{equation}

We use induction. Note that (\ref{37}) is valid for $x=2$.  Assuming that it is valid
for $x=2m,$ we have

$$
\sum_{0\leq n< 2m+2}(-1)^{s(n)}= (-1)^{s(2m)}+
(-1)^{2m+1}=0.
$$

Thus, if $x$ and $y$ are even, then, by (\ref{35})-(\ref{37}), we find

$$
3S_{3,0}([x,y))- S_{3,0}([4x, 4y))= S_{1,0}([x, y))=0.
$$

and

\begin{equation}\label{38}
S_{3,0}([4x, 4y))= 3S_{3,0}([x, y)).
\end{equation}
\newpage
Now, using (\ref{38}), we find

$$
S_{3,0}(2^m)= 3 S_{3,0}(2^{m-2})=3^2 S_{3,0}(2^{m-4})=\ldots=
$$
$$
=\begin{cases} 3^{\frac m 2 -1} S_{3,0}([0,4)), \;\; if \;\; m
\;\;is \;\; even\\3^{\frac{m-1}{2}}S_{3,0}([0, 2)),\;\; if \;\; m
\;\; is \;\; odd \end{cases}=
$$

$$
=\begin{cases} 2\cdot 3^{\frac m 2 -1},\;if \; m \; is \;even,\\
3^{\frac{m-1}{2}}\;if \; m \; is \; odd \end{cases}
$$

which proves (\ref{21}). Furthermore,

$$S_{3,0}([2^n, 2^n+2^m))=3S_{3,0}(2^{n-2},2^{n-2}+2^{m-2})=$$$$3^2S_{3,0}([2^{n-4},
2^{n-4}+2^{m-4}))= \ldots=
$$
$$
=\begin{cases} 3^{\frac m 2 -1}\cdot S_{3,0}([2^{n-m+2},
2^{n-m+2}+4)), \;\; if \;\; m \;\;is \;\;
even\\3^{\frac{m-1}{2}}S_{3,0}([2^{n-m+1}, 2^{n-m+1}+2)),\;\; if
\;\; m \;\; is \;\; odd.
\end{cases}
$$

Since

$$
2^a\equiv \begin{cases} -2,\;\; if \;\; a \;\; is \;\; even\\
-1,\;\; if \;\; a \;\; is \;\; odd \;\;(mod 3),\end {cases}
$$
we have

$$
S_{3,0}([2^a, 2^a+2))=\begin{cases} 0,\;\; if \;\;a \;\; is \;\;
even,\\
1, \;\; if \;\; a \;\; is \;\; odd, \end {cases}
$$

$$
S_{3,0}([2^a, 2^a + 4))= 1.
$$

This completes proof of the theorem. $\blacksquare$

\begin{remark}\label{r2}  Formulas (\ref{31}) and (\ref{34}) express $S_{3,i}(x),\;\; i=1,2$
 via $S_{3,0}(x)$. Thus, the above algorithm for calculation of
$S_{3,0}(x)$ allows to calculate also $S_{3,1}(x)$ and $S_{3,2}(x).$
\end{remark}

\begin{remark}\label{r3}  Our algorithm allows to calculate also the sums of
the form $S_{3\cdot 2^m,r}(x)$ . For example, if \; $0\leq r\leq
2^m-1,\;\; n> m,$  then we have

$$
S_{3\cdot 2^m, k\cdot 2^m+r}(2^n)=(-1)^{s(r)} S_{3\cdot 2^m,
k\cdot 2^m}(2^n)=$$ $$(-1)^{s(r)}S_{3,k}(2^{n-m}),\;\;k=0,1,2.
$$
\end{remark}
\newpage
\section{Lower estimate}

     For integers $x,y,$ such that $0\leq x< y,$ denote

\begin{equation}\label{39}
\delta_{3,0}([x,y))=\frac{S_{3,0}([x,y))}{(x-y)^\lambda}.
\end{equation}

\begin{lemma}\label{l1}

If $t(x)\equiv 0\;\pmod 6,$ then, for $0\leq x < y,$ we have

\begin{equation}\label{40}
\delta_{3,0}([x,y))\geq \delta_{3,0}([0,x))+\delta_{3,0}([0,y-x)).
\end{equation}
\end{lemma}
\slshape Proof.\upshape  Indeed, evidently, we have

\begin{equation}\label{41}
S_{3,0}([0,y))= S_{3,0}([0,x))+S_{3,0}([x,y)).
\end{equation}

Now since $t(x)\equiv 0\;(\mod 6),$ then, in view of the bijection of shift
$[x,y)\leftrightarrow [0,y-x),$ the numbers divisible by 3 in $[x,y)$
correspond to multiple of 3 in $[0, y-x)$ with the
conservation of parity of the binary sums. Therefore, by
(\ref{46}), we have

$$
S_{3,0}([0,y))=S_{3,0}([0,x))+S_{3,0}([0,y-x))
$$

and thus
\begin{equation}\label{42}
\delta_{3,0}([0,y))= \delta_{3,0}([0,x))\frac{x^\lambda}{y^\lambda}+
\delta_{3,0}([0,y-x))\frac{(y-x)^\lambda}{y^\lambda}.
\end{equation}

Using the Jensen inequality for the convex function $x^\lambda$, we
have

$$
x^\lambda+(y-x)^\lambda\geq y^\lambda
$$

and, consequently, (\ref{42}) yields (\ref{40}). $\blacksquare$

Denote by $\inf^*_{s(N)\geq 3} \delta_{3,0}([0,N))$ the infimum
of $\delta([0,N))$ over all $N$ with $s(N)\geq 3$ which have
either all the even exponents of 2 or all the odd exponents of 2 in
its binary expansion. Taking in Lemma \ref{l1} in the capacity of
$x$ the sums of the forms $2^n+2^m,\;n-m\equiv 1\pmod
2,\;\;\sum^6_{i=1}2^{2k_i},\;\sum^6_{i=1}2^{2k_i-1}$ we conclude
that

$$
\inf_{N\geq 2}\delta_{3,0}([0,N))\geq \min(\inf_{1\leq \sigma(N)\leq
2}\delta([0,N)),\;\inf^*_{3\leq \sigma(N)\leq 6}\delta([0,N)).
$$

Since the inverse inequality is evident, we have

\begin{equation}\label{43}
\inf_{N\geq 2}\delta_{3,0}([0,N))= \min(\inf_{1\leq \sigma(N)\leq
2}\delta([0,N)),\;\inf^*_{3\leq \sigma(N)\leq 6}\delta([0,N)).
\end{equation}
\newpage
Let first the binary expansion of $N$ have all the even (odd)
exponents of 2 with a fixed two largest exponents:

$$
N= 2^n + 2^m +\ldots,\;\;m<n,\;3\leq s(N)\leq 6
$$

According to Theorem \ref{t1}, formulas (\ref{15})-(\ref{20})
and (\ref{43}), we easily find

$$
S_{3,0}(N)\geq 2\cdot 3^{\frac n 2 -1}+3^{\frac m 2
-1}-3^{\frac{m-2}{2}-1}-2\cdot
3^{\frac{m-4}{2}-1}-3^{\frac{m-6}{2}-1}>
$$
\begin{equation}\label{44}
>2\cdot 3^{\frac n 2 -1},\; if \; n\;and\; m\; are\;even,
\end{equation}

$$
S_{3,0}(N)\geq 3^{\frac {n-1}{2}}+3^{\frac
{m-1}{2}}-0-3^{\frac{m-5}{2}} -3^{\frac{m-7}{2}}>
3^{\frac{n-1}{2}}-3^{\frac{m-1}{2}}+
$$

\begin{equation}\label{45}
+3^{\frac {m-1}{2
}}(\frac{2}{\sqrt{3}}-\frac{1}{9\sqrt{3}}-\frac{1}{27\sqrt{3}}),\;
if \; n\;and\; m\; are\;odd.
\end{equation}

In the case of (\ref{44})

$$
\delta_{3,0}{N}>\frac{2\cdot 3^{\frac n
2-1}}{(2^n+2^m+2^{m-2}+2^{m-4}+\ldots)^\lambda}>
$$
\begin{equation}\label{46}
>\frac{\frac 2 3\cdot 3^{\frac n
2}}{(2^n+\frac{2^{m+2}}{3})^\lambda}>\frac 2 3 \frac{3^{\frac n
2}}{2^{n\alpha}}=\frac 2 3,
\end{equation}

while in the case of (\ref{45})

$$
\delta_{3,0}(N)>\frac{\left(3^{\frac{n-1}{2}}-3^{\frac{m-1}{2}}\right)+3^{\frac{m-1}{2}}\cdot
1.069...}{\left(2^n+\frac{2^{m+2}}{3}\right)^\lambda}>
$$
$$
>\frac{3^{\frac{n-1}{2}}-3^{\frac{m-1}{2}}}{2^{\eta\lambda}}\cdot
\frac{2^{\eta\lambda}}{\left(2^n+\frac{2^{m+2}}{3}\right)^\lambda}+\frac{3^{\frac{m-1}{2}}1.069...}
{\left(\frac{2^{m+2}}{3}\right)^\lambda}\cdot
\frac{\left(\frac{2^{m+2}}{3}\right)^\lambda}{\left(2^n+\frac{2^{m+2}}{3}\right)^\lambda},
$$

and, using the Jensen inequality, we have

$$
\delta_{3,0}(N)>\min\left(\left(1-\frac{1}{3^{\frac{n-m}{2}}}\right),\;
1.069...\;\frac{3^{\frac{m-1}{2}+\lambda}}{3^{\frac m 2 +1}}
\right)\geq
$$
\begin{equation}\label{47}
\geq \min\left( \frac 2 3,\;\frac{1.069...}{3^{1.5-\lambda}}\right)>
0.49.
\end{equation}

Now we show that

\begin{equation}\label{48}
\inf_{1\leq\sigma(N)\leq2,N\geq4}\delta_{3,0}([0,N))=
\frac{2}{6^\lambda}= 0.48345...\;.
\end{equation}

We consider several cases.
\newpage
a) $s(N)=1$. Taking into account that $4^\lambda=3,$ note that,
for even $x,y,$ (\ref{38}) yields

\begin{equation}\label{49}
\delta_{3,0}([4x, 4y))= \delta_{3,0}([x, y)).
\end{equation}

Consequently, for $n\geq 1,$ we find

\begin{equation}\label{50}
\delta_{3,0}([0,2^n))= \begin{cases} \delta_{3,0}([0,4),\; if\; n\;
is\;even\\\delta_{3,0}([0,2),\; if\; n\; is\;odd \end{cases}=
\begin{cases} \frac 2 3,\; n\; is\; even\\\frac{1}{\sqrt{3}}, \; n
\; is\; odd \end{cases}
\end{equation}

$b)$ $s(N)=2$.  Note that, in the cases $N=2^n+2^m,$ where $n$ and
$m$ are of the same parity, according to (\ref{44}) and (\ref{45}),
in any case, (\ref{47}) satisfies.

$b_1)$ Let $n$ be odd, $m$ be even, $n> m$. Then

\begin{equation}\label{51}
S_{3,0}([0,2^n+2^m))=S_{3,0}([0,2^n))+S_{3,0}([2^n, 2^n+2^m))
\end{equation}

If here $m=0,\; n\geq 3,$ then, according to Theorem 1,

\begin{equation}\label{52}
S_{3,0}([0,2^n+1))=S_{3,0}([0,2^n))=3^{\frac{n-1}{2}}
\end{equation}
and
$$
\delta_{3,0}([0,2^n+1 ))=\frac{3^\frac{n-1}{2}}{(2^n+1)^\lambda}=$$ $$
\frac{3^\frac{n-1}{2}}{3^{\frac n 2}(1+\frac{1}{2^n})^\lambda}\geq
\frac{1}{\sqrt{3}(1+\frac 1 8)^\lambda}= 0.52589...\;.
$$

Note that, the case $m=0, n=1$ corresponds to $N=3,$ and now we do not
consider this case (see (\ref{48})).

     Let now in (\ref{51}) $m\geq 2,\;n\geq 3$. Then, by Theorem \ref{t1}, we have

$$
S_{3,0}([0, 2^n+2^m))= 3^\frac{n-1}{2} + 3^ \frac{m-2}{2}
$$

and
\begin{equation}\label{53}
\delta_{3,0}([0,
2^n+2^m))=\frac{3^\frac{m-2}{2}(3^\frac{n-m+1}{2}+1)}{2^{\lambda
m}(2^{n-m}+1)^\lambda}=\frac 1 3 \cdot \frac
{3^\frac{n-m+1}{2}+1}{(2^{n-m}+1)^\lambda}.
\end{equation}

Put $n-m=x\geq 1$. The function $f(x)=\frac 1 3
\frac{3^\frac{x+1}{2}+1}{(2^x +1)^\lambda}$  tends to
$\frac{1}{\sqrt{3}}=0.57735...$ for  $x\rightarrow\infty$ and
has the unique extremum in the point $x=\frac{\ln 3}{\ln 4 - \ln 3}=
3.8188... $  Besides, $f(1)= 0.5582...,\;f(3) = 0.5843..., \; f(5)= 0.5843...\;.$ Thus, for odd $n$  and  even $m,$ we
have
\begin{equation}\label{54}
\delta_{3,0}([0, 2^n+2^m))\geq 0.52589...\;.
\end{equation}
\newpage
$b_2)$ Let $n$ be even, $m$ be odd, $n>m\geq 1.$ Then, according to
(\ref{51}) and Theorem \ref{t1}, we have

$$
S_{3,0}([0, 2^n + 2^m))= 2\cdot 3^{\frac n 2 -1}
$$

and

\begin{equation}\label{55}
\delta_{3,0}([0, 2^n+2^m))=\frac{2\cdot 3^\frac n 2-1}{2^{\lambda
m}(2^{n-m}+1)^\lambda}=\frac 2 3 \cdot \frac
{3^\frac{n-m}{2}}{(2^{n-m}+1)^\lambda}.
\end{equation}

Since the function $g(x)=\frac 2 3\frac{3^\frac x
2}{(2^x+1)^\lambda}$ is increasing and tends to $\frac 2 3,$ while
$g(1)=\frac{2}{\sqrt{3}\cdot 3^\lambda}=\frac{2}{6^\lambda}=
0.48345...,$ then (\ref{48}) is proved. Moreover, the infimum
(\ref{53}) is realized on $\delta_{3,0}([0, 2^n + 2^{n-1})),\;n=2,4,6,
\ldots\;.$

Finally, note that

$$
\delta_{3,0}([0,1))=1,
\;\;\delta_{3,0}([0,2))=\frac{1}{2^\lambda}=\frac {1}{\sqrt{3}}=
0.577...,
$$
\begin{equation}\label{56}
\delta_{3,0}([0,3))=\frac {1}{\sqrt{3}}=0.41868...\;.
\end{equation}

From (\ref{48}) and (\ref{56}) it follows that

\begin{equation}\label{57}
\inf_{N\geq 1}\delta_{3,0}([0,N))=\frac {1}{3^\lambda}=0.41868...
\end{equation}

and
\begin{equation}\label{58}
\liminf_{N\rightarrow\infty}\delta_{3,0}([0,N))=\frac
{2}{3^\lambda}=0.48345...,
\end{equation}

such that $\liminf_{N\rightarrow\infty}\delta_{3,0}([0,N))$ is
realized on the sequence

\begin{equation}\label{59}
N_n= 2^n + 2^{n-1},\;\;n=2,4,6 \ldots\;.
\end{equation}

Finally, we obtain the exact lower bounds for $S_{3,0}(N).$ For
$N\geq 4,$ we have

$$
S_{3,0}(N)\geq 2\left(\frac N 6\right)^\lambda=0.48345...
N^\lambda.
$$
In particular, for $m\geq 2,$ we obtain:
$$S_{3,0}(3n)\geq 2\left(\frac n 2\right)^\lambda=\frac
{2}{\sqrt{3}}n^\lambda=1.1547005... n^\lambda
$$
that corresponds to (\ref{11}).
In addition, note that, for all $N\geq 1,$
\begin{equation}\label{60}
S_{3,0}(N)\geq \left\lfloor 2\left(\frac N 6\right)^\lambda
\right\rfloor.
\end{equation}
\newpage

\section{Upper estimate}

Let the binary expansion of even $N$ is

\begin{equation}\label{61}
N=2^{n_1}+2^{n_2}+ \ldots +2^{n_r},\;\; n_1> n_2 >\ldots >n_r\geq 1.
\end{equation}

     Using Theorem \ref{t1} and formulas (\ref{15})-(\ref{20}), we easily
find

\begin{equation}\label{62}
S_{3,0}(N)\leq S_{3,0}(2^{n_1}+2^{n_2}),\;\; if \;\;2\leq
s(N)\leq 5.
\end{equation}

Let now $s(N)\geq 6$. Then we have

$$
S_{3,0}(N)\leq S_{3,0}(2^{n_1}+2^{n_2})-0-2\cdot
3^{\frac{n_4}{2}-1}-0+2\cdot 3^{\frac{n_6}{2}-1}+
$$

$$
+3^\frac{n_7-1}{2}+2\cdot 3^{\frac{n_8}{2}-1}+ 3^\frac{n_9-1}{2}+
2\cdot 3^{\frac{n_{10}}{2}-1}.
$$

Note that $n_6\leq n_4-2,\;n_7\leq n_4-3,\ldots\;.$  Hence,

$$
S_{3,0}(N)\leq S_{3,0}(2^{n_1}+2^{n_2})-0-2\cdot
3^{\frac{n_4}{2}-1}+2\cdot 3^{\frac{n_4-2}{2}-1}+3^\frac{n_4-4}{2}+
$$

$$
+2\cdot 3^{\frac{n_4-4}{2}-1}+3^\frac{n_4-6}{2}+2\cdot
3^{\frac{n_4-6}{2}-1}+\ldots=S_{3,0}(2^{n_1}+2^{n_2})-
$$

$$
-2\cdot3^{\frac{n_4}{2}-1}+\frac 2
3\cdot\frac{3^{\frac{n_4}{2}-1}}{\frac 2 3}+\frac 1 9\cdot
\frac{3^{n_4}}{\frac 2 3}=
$$

\begin{equation}\label{63}
=S_{3,0}(2^{n_1}+2^{n_2})-3^{\frac{n_4}{2}}(\frac 2 3-\frac 1
3-\frac 1 6)<S_{3,0}(2^{n_1}+2^{n_2}).
\end{equation}
Thus, if $N$ is even, then
\begin{equation}\label{64}
\delta_{3,0}(N)<\delta_{3,0}(2^{n_1}+2^{n_2}).
\end{equation}
Let us show that (\ref{64}) is correct also for odd $N$ with
$s(N-1)\geq 2$ except for, probably, the case of
$s(N-1)=2,\;\; N\equiv 1 \pmod 3.$ Indeed, if $s(N-1)=2,$ then
it follows directly from (\ref{23}). If $s(N-1)=3,$ such that
$N=2^{n_1}+2^{n_2}+2^{n_3}+1,$ then, by (\ref{23}), also $S_{3,0}(N)\leq
S_{3,0}(N-1),$ and, by (\ref{64}), for $N-1,$ we have
$$
\delta_{3,0}(N)\leq \delta_{3,0}(N-1)<\delta_{3,0}(2^{n_1}+2^{n_2}).
$$
Finally, if $s(N-1)=4$ and $s(N-1)\geq 5,$ then the
inequalities (\ref{62}) and (\ref{63})
 again give (\ref{64}) also for odd $N.$\newline
\indent Let us find $\sup_{\sigma(N)\leq 2}\delta(N)$. By (\ref{50}),
$\sup_{\sigma(N)= 1,N\geq 2}\delta(N)=\frac 2 3.$
Note that
\begin{equation}\label{65}
\delta_{3,0}(1)=1.
\end{equation}
\newpage
Moreover, we have already investigated $\delta(2^n+2^m),$ if
$n-m\equiv 1 \pmod 2$ (see(\ref{53}) and  (\ref{55})) and know that
in this case $\delta(2^n+2^m)\leq \frac 2 3.$ Let us consider the
remaining cases.
a) Let $n$ and $m$ be odd, $1\leq m\leq n-2.$ Then, by Theorem
\ref{t1},
$$
S_{3,0}(2^n+2^m)=3^{\frac{n-1}{2}}+3^{\frac{m-2}{2}}
$$
and thus
$$
\delta_{3,0}(2^n+2^m)=\frac{3^\frac{n-1}{2}+3^\frac{m-2}{2}}{(2^n+2^m)^\lambda}=
\frac 1 3\cdot \frac{3^{\frac{n-m+1}{2}}+1}{(2^{n-m}+1)^\lambda}.
$$
This coincides with (\ref{53}) for $n-m$ is even. Thus
$$
\delta_{3,0}(2^n+2^m)\leq \frac 1 3
\cdot\frac{3^{1.5}+1}{5^\lambda}=0.5768... <\frac 2 3.
$$
b) Let $n$ and $m$ be even, $0\leq m \leq n-2$. In the case $m=0,$ by
(\ref{23}) and (\ref{50}), we have $\delta_{3,0}(2^n+1)< \frac 2 3.$
Let now $m\geq 2$. Then, by Theorem \ref{t1},
$$
\delta_{3,0}(2^n+2^m)=\frac{2\cdot 3^{\frac n 2 -1}+3^{\frac m 2
-1}}{(2^n+2^m)^\lambda}= \frac 1 3\cdot
\frac{2\cdot3^{\frac{n-m}{2}}+1}{(2^{n-m}+1)^\lambda}.
$$
Put $n-m=x\geq 2$. The function $h(x)=\frac 1 3\cdot \frac{2\cdot
3^\frac x 2+1}{(2^x+1)^\lambda}$ tends to $\frac 2 3$ as
$x\rightarrow\infty$ and has the unique extremum in the point
$x=\frac{\ln 2}{\ln 2-\ln\sqrt{3}}=4.8188...$  We have
$$h(2)=0.6517...,\; h(4)=0.67069...,\;h(6)=\frac 1 3\cdot
\frac{55}{{65}^\lambda}=0.67072... > \frac 2 3.$$
Thus, according to (\ref{65}),
$$
\sup_{N\geq 1}\delta_{3,0}(N)=1,
$$
if $N=1$ is not isolated. If $N\geq 2$ (the case when $N$ is odd and
$s(N)=3$ we consider in the supposition $N\equiv 0\; or\; 2\;\pmod 3$),
the supremum is $\frac 1 3\frac {55}{{65}^\lambda}=0.67072...$
In this case
$$
S_{3,0}(N)\leq \lfloor 0.67072... N^\lambda\rfloor.
$$
In particular,
$$
S_{3,0}(3N)\leq \frac{55}{3}(\frac{3}{65})^\lambda N^\lambda
$$
that corresponds to (\ref{11}).\newline
\indent Finally, in the case $N\equiv 1\pmod 6$ with $s(N)=3,$ we
have
$$
S_{3,0}(N)\leq S_{3,0}(N-1)+1\leq $$ $$\lfloor
0.67072...(N-1)^\lambda\rfloor+1\leq\lceil 0.67072...
N^\lambda\rceil.
$$
This equality holds, for example, on $N=19$ and $67.$
Now we conclude that
$$
\limsup_{N\rightarrow\infty}\delta_{3,0}(N)=\frac 1 3 \cdot\frac
{55}{(65)^\lambda}
$$
and it is realized on the sequence $2^n+2^{n-6}$ for even $n.$\newline
\indent Thus, we proved that, for $N\geq 2,$
\newpage
$$
\lfloor 0.483459... N^\lambda\rfloor \leq S_{3,0}(N)\leq\lceil
0.67072... N^\lambda\rceil.
$$

\section{Fast computing algorithm}

Now we give a fast algorithm for evaluation of
$S_{3,0}(N)$ which is based on the following formulas:
\begin{equation}\label{66}
S_{3,0}(N)=3S_{3,0}\left(\left\lfloor \frac N
4\right\rfloor\right)+\nu(N),
\end{equation}

where

\begin{equation}\label{67}
\nu(N)=\begin{cases} 0,\;if \;N\equiv 0,7,8,9,16,17,18,22,23\;\pmod
{24};\\ (-1)^{s(N)},\; if \; N\equiv 3,4,10,12,20 \; \pmod {24};\\
(-1)^{s(N)+1},\; if \; N\equiv 1,2,5,6,11,19,21\; \pmod {24};\\
2(-1)^{s(N)},\; if \; N\equiv 15\;\pmod {24};\\
2(-1)^{s(N)+1},\; if \; N\equiv 13,14\;\pmod {24}.\end{cases}
\end{equation}

Note that, by the definition (see (\ref{3})), $S_{3,0}(0)=0$. This
algorithm is very suitable for sharp experiments. For example, let us show the concluding phase of the
calculation of  $S_{3,0}(500000)$ by this algorithm:

$$
S_{3,0}(500  000)= 6561
S_{3,0}(7)-2187(-1)^{s(30)}-
$$

$$
-729(-1)^{s(122)}+27(-1)^{s(7812)} -9(-1)^{s(31250)}=
$$

$$
=19683-2187+729+27+9=18261.
$$
\section{Proof of formulas (66)-(67)}
We use induction. Evidently, formulas (\ref{66})-(\ref{67}) are valid for $x=1.$
Suppose that
$$S_{3,0}(x)=3S_{3,0}\left(\left\lfloor \frac x
4\right\rfloor\right)+\nu(x)$$
is true and prove that
$$S_{3,0}(x+1)=3S_{3,0}\left(\left\lfloor \frac {x+1}
4\right\rfloor\right)+\nu(x+1).$$
Note that
$$S_{3,0}(x+1)-S_{3,0}(x)=\begin{cases} 0,\;\; if \; x \;\; is\;\; not\;\;
multiple\;\; of\;\;3,\\
(-1)^{s(x)}, \;\; if \;\; x \;\; is\;\; multiple \;\; of\;\; 3.\end{cases}$$
Therefore, it is sufficient to prove that
$$3(S_{3,0}\left(\left\lfloor \frac {x+1}
4\right\rfloor\right)-S_{3,0}\left(\left\lfloor \frac {x}
4\right\rfloor\right))+\nu(x+1)-\nu(x)=$$
\newpage
\begin{equation}\label{68}
\begin{cases} 0,\;\; if \; x \;\; is\;\; not\;\;
multiple\;\; of\;\;3,\\
(-1)^{s(x)}, \;\; if \;\; x \;\; is\;\; multiple \;\; of\;\; 3.\end{cases}
\end{equation}
\begin{lemma}\label{l2}
We have
$$S_{3,0}(\lfloor \frac {x+1}
4\rfloor)-S_{3,0}(\lfloor \frac {x}
4\rfloor)=$$
\begin{equation}\label{69}
\begin{cases} (-1)^{s(3t)},\;\; if \; x \;\; has\;\; the\;\;
form\;\; x=12t+3,\\
0, \;\; otherwise.\end{cases}
\end{equation}
\end{lemma}
\slshape Proof.\;\;\upshape Let $x\equiv j \pmod{4}.$ Then
$$\lfloor \frac {x+1}4\rfloor-\lfloor \frac {x}
4\rfloor=\begin{cases} 1,\;\; if \; j=3,
\\0, \;\; otherwise.\end{cases} $$
Thus
$$S_{3,0}(\lfloor \frac {x+1}
4\rfloor)-S_{3,0}(\lfloor \frac {x}
4\rfloor)=$$
$$\begin{cases} (-1)^{s(\lfloor \frac {x}
4\rfloor)},\;\; if \; j=3\;\; and\;\;\lfloor \frac {x}
4\rfloor\equiv0\pmod3, \\
0, \;\; otherwise,\end{cases}$$
and the lemma follows. $\blacksquare$\newline
Hence, by (\ref{68}) and Lemma \ref{l2}, it is left to prove that
$$\nu(x+1)-\nu(x)=$$
\begin{equation}\label{70}
\begin{cases} (-1)^{s(x)},\;\; if \; x \;\; has\;\; the\;\;
form\;\; x=12t, 12t+6\; or\; 12t+9,\\
(-1)^{s(x)}-3(-1)^{s(3t)}, \;\; if \; x \;\; has\;\; the\;\; form\;\;
 x=12t+3,\\0, \;\; if \; x \;\; is\;\; not\;\; multiple\;\; of \;\;3.\end{cases}
\end{equation}
Distinguish the corresponding cases.
 By the definition of $\nu(x),$ we should consider the cases
 $$a)\;\; x=0,6,9,12,18,21;$$  $$b)\;\;  x=3,15$$ and
  $$c)\;\; x=1,2,4,5,7,8,10,11,13,14,16,17,19,20,22,23.$$
  Below we use the definition (\ref{67}) of $\nu(x)$ for the check of (\ref{70}).\newline
$$Case \;\;a)$$  $$x=0,\;\nu(1)-\nu(0)=(-1)^{s(1)+1}-0=(-1)^{s(0)};$$
  $$x=6,\;\nu(7)-\nu(6)=0-(-1)^{s(6)+1}=(-1)^{s(6)};$$
  $$x=9,\;\nu(10)-\nu(9)=(-1)^{s(10)}-0=(-1)^{s(9)};$$
  $$x=12,\;\nu(13)-\nu(12)=2(-1)^{s(13)+1}-(-1)^{s(12)}=(-1)^{s(12)};$$
  $$x=18,\;\nu(19)-\nu(18)=(-1)^{s(19)+1}-0=(-1)^{s(18)};$$
  \newpage
  $$x=21,\;\nu(22)-\nu(21)=0-(-1)^{s(21)+1}=(-1)^{s(21)}.$$
$$Case \;\;b)$$  $$x=3\;(t=0),\;\nu(4)-\nu(3)=-2=(-1)^{s(3)}-3(-1)^{s(0)};$$
  $$x=15\;(t=1),\;\nu(16)-\nu(15)=0-2(-1)^{s(15)}=(-1)^{s(15)}-3(-1)^{s(3)}.$$

  $$Case \;\;c)$$
  $$x=1,\;\nu(2)-\nu(1)=(-1)^{s(2)+1}-(-1)^{s(1)+1}=0;$$
  $$x=2,\;\nu(3)-\nu(2)=(-1)^{s(3)}-(-1)^{s(2)+1}=0;$$
  $$x=4,\;\nu(5)-\nu(4)=2(-1)^{s(5)+1}-(-1)^{s(4)}=0;$$
  $$x=5,\;\nu(6)-\nu(5)=(-1)^{s(6)+1}-(-1)^{s(5)+1}=0;$$
  $$x=7,\;\nu(8)-\nu(7)=0-0=0;$$
  $$x=8,\;\nu(9)-\nu(8)=0-0=0;$$
  $$x=10,\;\nu(11)-\nu(10)=(-1)^{s(11)+1}-(-1)^{s(10)}=0;$$
  $$x=11,\;\nu(12)-\nu(11)=(-1)^{s(12)}-(-1)^{s(11)+1}=0;$$
  $$x=13,\;\nu(14)-\nu(13)=2(-1)^{s(14)+1}-2(-1)^{s(13)+1}=0;$$
  $$x=14,\;\nu(15)-\nu(14)=2(-1)^{s(15)}-(-1)^{s(14)+1}=0;$$
  $$x=16,\;\nu(17)-\nu(16)=0-0=0;$$
  $$x=17,\;\nu(18)-\nu(17)=0-0=0;$$
  $$x=19,\;\nu(20)-\nu(19)=(-1)^{s(20)}-(-1)^{s(19)+1}=0;$$
  $$x=20,\;\nu(21)-\nu(20)=(-1)^{s(21)+1}-(-1)^{s(20)}=0;$$
  $$x=22,\;\nu(23)-\nu(22)=0-0=0;$$
   $$x=23,\;\nu(24)-\nu(23)=0-0=0.\; \blacksquare$$
   Now from (\ref{66})-(\ref{67}) we obtain
\begin{corollary}\label{col71}The function
\begin{equation}\label{71}
F(x)=(-1)^{s(x)}(S_{3,0}(x)-3S_{3,0}(\lfloor \frac x
4\rfloor))
\end{equation}

is periodic with the period $24.$
\end{corollary}
\newpage
\section{An expression for $\nu(x)$}
Let us back to (\ref{35})-(\ref{37}). For $x=0, y:=x,$ we have
$$3S_{3,0}(x)-S_{3,0}(4x)=$$
\begin{equation}\label{72}
\sum_{0\leq i< x}(-1)^{s(i)}=\begin{cases} 0,\;\; if
\; x \;\; is\;\; even,\\
(-1)^{s(x-1)}, \;\; if \; x\;\; is\;\;odd.\end{cases}
\end{equation}
From this we have
\begin{corollary}\label{col2} The function
\begin{equation}\label{73}
G(x)=(-1)^{s(x-1)}(3S_{3,0}(x)-S_{3,0}(4x))
\end{equation}
is periodic with the period $2.$
\end{corollary}
Write (\ref{73}) for
$x:=\lfloor\frac{x}{4}\rfloor.$ We have
$$(-1)^{s(\lfloor\frac{x}{4}\rfloor-1)}(3S_{3,\;0}(\lfloor\frac{x}{4}\rfloor)-
S_{3,\;0}(4\lfloor\frac{x}{4}\rfloor)= $$
\begin{equation}\label{74}
\begin{cases} 0,\;\; if \; \lfloor\frac{x}{4}\rfloor \;\; is\;\; even,\\
1, \;\; if \;\; \lfloor\frac{x}{4}\rfloor \;\; is\;\; odd.\end{cases}
\end{equation}
Note that $\lfloor\frac{x}{4}\rfloor$ is even, if
$x=0,1,2,3,8,9,10,11,...$ and odd for other integers. Thus we obtain
\begin{corollary}\label{col3} The function
\begin{equation}\label{75}
H(x)=(-1)^{s(\lfloor \frac x
4\rfloor-1)}(3S_{3,\;0}(\lfloor\frac{x}{4}\rfloor)-
S_{3,\;0}(4\lfloor\frac{x}{4}\rfloor)
\end{equation}
is periodic with the period $16,$ such that
\begin{equation}\label{76}
 H(x)=\begin{cases} 0,\;\; if \; x\equiv0,1,2,3,8,9,10,11\pmod{16},\\
1, \;\; if \;\;x\equiv4,5,6,7,12,13,14,15\pmod{16}.\end{cases}
\end{equation}
\end{corollary}
Consider the difference
\begin{equation}\label{77}
\delta(x)=S_{3,0}(x)-S_{3,0}(4\lfloor \frac {x}
4\rfloor).
\end{equation}
\begin{lemma}\label{l3}
We have
\begin{equation}\label{78}
\delta(x)=\begin{cases} (-1)^{s(x-1)},\;\; if \; x\equiv1,7\;or\;10\pmod{12}\\
 (-1)^{s(x-2)}, \;\; if \; x\equiv2\;or\;11\pmod{12} \\(-1)^{s(x-3)}, \;\; if
 \; x\equiv3 \pmod{12}\\0, otherwise. \end{cases}
 \end{equation}
\end{lemma}
\slshape Proof.\;\;\upshape Let $x=12t+j,\;j=0,1,...,11.$ Consider 3 cases.
$$a)\; j=0,1,2\;or\;3.$$
Then
$$\delta(x)=S_{3,0}(12t+j)-S_{3,0}(12t)=$$
\newpage
 $$\begin{cases} 0,\;\; if \; j=0,\\
(-1)^{s(x-j)}, \;\; if \;\; j=1,2,3.\end{cases}$$

$$b)\; j=4,5,6\;or\;7.$$
Then
 $$\delta(x)=S_{3,0}(12t+j)-S_{3,0}(12t+4)=$$

 $$\begin{cases} 0,\;\; if \; j=4,5,6,\\
(-1)^{s(x-1)}, \;\; if \;\; j=7.\end{cases}$$

$$c)\; j=8,9,10\;or\;11.$$
Then
 $$\delta(x)=S_{3,0}(12t+j)-S_{3,0}(12t+8)=$$

  $$\begin{cases} 0,\;\; if \; j=8,9,\\
(-1)^{s(x-1)}, \;\; if \;\; j=10,\\(-1)^{s(x-2)}, \;\; if \;\; j=11\end{cases}\;
$$
and (\ref{78}) follows.\; $\blacksquare$\newline
 Using (\ref{75}) and (\ref{77}), we obtain the following expression for $\nu(x).$
 \begin{corollary}\label{col4}
 \begin{equation}\label{79}
\nu(x)=\delta(x)-(-1)^{s(\lfloor \frac {x}
4\rfloor-1)}H(x),
\end{equation}
where $H(x)$ and $\delta(x)$ are defined by $(\ref{76})$ and $(\ref{78})$
 respectively.
 \end{corollary}
Finally note that, using (\ref{79}) and some simple transformations, one can obtain
the second proof of formulas (\ref{66})-(\ref{67}).

\;\;\;\;\;\;

\end{document}